\documentclass[12pt,a4paper]{article}
\usepackage{graphicx}
\usepackage{latexsym}
\usepackage{amsmath}
\usepackage{amssymb}
\usepackage{mathrsfs}
\usepackage{array} 

\numberwithin{equation}{section}

 \date{}
\begin{document}
 \title{ Local geometric proof of Riemann Hypothesis}
  \author{\large  CHUANMIAO CHEN \\
    {\small School of Mathematics and Statistics, Central South University, Changsha, China,}\\
    {\small Email: cmchen@hunnu.edu.cn}  }

   \maketitle
 {\small
    \begin{center}
   \begin{minipage}{5.8in}
    {\bf Abstract}\ \
       Riemann function $\xi(s)=u+iv, s=\beta+1/2+it$ has the important symmetry: $v=0$ if $\beta=0$.
    For $\beta>0$ we prove $|u|>0$ inside any root-interval $I_j=[t_j,t_{j+1}]$ and
    $v$ has opposite signs at two end-points of $I_j$.
    They imply local peak-valley structure and $||\xi||=|u|+|v/\beta|>0$ in $I_j$.
    Because each $t$ must lie in some $I_j$, then $||\xi||>0$ is valid for any $t$.
    By the equivalence $Re(\frac{\xi'}{\xi})>0$ of Lagarias(1999),
   we show that RH implies the peak-valley structure,
   which may be the geometric model expected by Bombieri(2000).
   {\bf AMS  Classification of Subjects}  11M26, 65E05\\
  \end{minipage}
 \end{center}
   }

\section{Introduction. Difficulty and hope}

 In 1737 Euler proved that the product formula of the prime number $p$
    \begin{equation}
   \zeta(s)=\sum_{n=1}^{\infty}\frac{1}{n^s}=\prod_{p\ prime}(1-\frac{1}{p^s})^{-1}
   \end{equation}
   is convergent for $Re(s)>1$, but divergent for $Re(s)\le 1$. In 1859
  Riemann considered the complex variable $s=\sigma+it,\sigma>1$, using Gamma
   function $\Gamma(s/2)$, and got
   \[
   \zeta(s)=\sum_{n=1}^{\infty}n^{-s}
   =\pi^{s/2}\Gamma^{-1}(\frac{s}{2})\int_0^{\infty}x^{s/2-1}\psi(x)dx,\
   \ \psi(x)=\sum\limits_{n=1}^{\infty}e^{-n^2\pi x}.
    \]
   Using the equality of Jacobi function $\psi(x)$
    \begin{equation}
    2\psi(x)+1=x^{-1/2}(2\psi(\frac{1}{x})+1),
    \end{equation}
    taking $z=1/x$ and transforming the integral
    \[
      \int_0^1z^{s/2-1}\psi(z)dz=\frac{1}{s(s-1)}
   +\int_1^{\infty}x^{-s/2-1/2}\psi(x)dx,
    \]
 Riemann derived {\bf the first expression}
   \begin{equation}
     \zeta(s)=\pi^{s/2}\Gamma^{-1}(\frac{s}{2})\{\frac{1}{s(s-1)}
   +\int_1^{\infty}(x^{s/2-1}+x^{-s/2-1/2})\psi(x)dx\},
   \end{equation}
   which is already analytically extended to the whole complex plane around $s=1$.
    \par
   Furthermore Riemann introduced the entire function
   \begin{equation}
   \xi(s)=\frac{1}{2}s(s-1)\pi^{-s/2}\Gamma(\frac{s}{2})\zeta(s),\ \ \xi(s)=\xi(1-s).
   \end{equation}
   Through replacing by $\zeta$ and integrating by parts twice, it follows that
   \begin{equation}
   \begin{array}{ccc}
   \xi(s)&= & \displaystyle \frac{1}{2}+\frac{s(s-1)}{2}\int_1^{\infty}(x^{s/2-1}+x^{-s/2-1/2})\psi(x)dx,\\
         & = & \displaystyle r_1+\int_1^{\infty} (x^{s/2-1}+x^{-s/2-1/2})(2x^2\psi''+3x\psi')dx,\\
  \end{array}
   \end{equation}
   where $r_1=\frac{1}{2}+\psi(1)+\frac{1}{4}\psi'(1)=0$ by (1.2).
   Riemann derived {\bf the second expression}
    \begin{equation}
    \label{xixi}
   \xi(s)=\int_1^{\infty}(x^{s/2-1}+x^{-s/2-1/2})f(x)dx,\ \ f(x)=2x^2\psi''+3x\psi',\\
   \end{equation}
   which is symmetric with respect to $s=1/2$. If $\sigma=1/2$, then $Im(\xi)=0$.
   \par
    Riemann thought that a number of zeros of $\zeta(s)$ in the critical region
    $\Omega=\{s=\sigma+it:\ 0\le \sigma \le 1, 0\le t<\infty\}$ has an estimate
   \begin{equation}
    N(T)=\frac{1}{2\pi}(T\ln\frac{T}{2\pi}-T)+O(\ln T),
   \end{equation}
   which was later proved by Mangoldt\cite{Borwein06} in 1905, then proposed the following hypothesis.
   \par
   {\bf  Riemann Hypothesis} (RH). {\em In the critical region
   $\Omega=\{s=\sigma+it:\ 0\le \sigma \le 1, 0\le t<\infty\}$,
    all the zeros of $\zeta(s)$ lie on the critical line $\sigma=1/2$, which is called
    the non-trivial zeros}.
   \par
   RH is an extremely difficult problem, which has stimulated the untiring research
   in the areas of the analytic number theory and the complex functions, even the scientific computation.
    Smale\cite{Smale00} in 1998 reported 18 mathematical problems
   for next century, which included RH.
   \par
    There have been many theoretical researches for RH\cite{Bombieri,Borwein06,Conrey03,Edwards}.
    A lot of numerical experiments verified that RH is valid. However RH has not been proved to be valid or false in theory.
   \par
   We can see from (1.7) that the average spacing between two zeros is less than
   $2\pi/\ln T$. To study the distribution of these zeros,
   there were lots of large scale numerical experiments, e.g. Lune et al. in \cite{Lune83,Lune86}
   searched out $1.5e+9$
   roots on the critical line where all roots were single, no double.
    These computations were finished by Euler-Maclaurin formula outside the critical line
    and Riemann-Siegel formula on the critical line. Here note that Riemann formula (1.6) has not been used.
   They emphasized that no nontrivial zeros were found in the critical
   strip $\{0\le\sigma\le 1,0\le t\le 5.6e+8\}$, which make people have the reason
   to believe RH is true.
   The authors listed lots of computed data and drew many curve figures,
   which have greatly inspired us to understand the function $\zeta(s)$.
   There have been two surprising phenomena on the critical line.
   \par
    1). There is a high peak in each subinterval of the curve, and $1\sim 9$ smaller peaks
   between two high peaks. They found the ratio of the high peak and low peak can reach 1000 times.
   \par
   2). There are $1\sim 8$ roots between two high peaks.  They found two pairs of large roots,
   in which two adjacent roots were very close to each other, and looked like a double root.
   \par
   It is likely that these terrible micro-structures have stopped the proof of RH
   by the pure analytical methods, but which have inspired us to consider
   the local geometry property of $\xi$.
    \par
   {\bf The difficulties and hope}.
   \par
   So far most studies have been focused on $\zeta$. There are the estimates in $t$ \cite{Edwards}(p.185,200)
   \begin{equation}
    \begin{array}{llll}
    |\zeta(\sigma+it)|&\le Ct^{1/4-\beta/2}\ln t,
    & 0\le \sigma\le 1,\ \ \beta=\sigma-1/2,\\
    |\zeta(\frac{1}{2}+it)|&=O(t^{\lambda}),& \lambda=1/6\ \ or\ \ \lambda=19/116,\\
    \end{array}
    \end{equation}
    which are possibly expressed as
    \begin{equation}
    \begin{array}{llll}
    |\zeta(\sigma+it)|&\le Ct^{1/6-\beta/3}\ln t,\ \ 0\le \beta=\sigma-1/2\le 1/2.\\
    \end{array}
    \end{equation}
   People can make more refined estimates, but it is not suitable to study its zeros.
     \par
    On the other hand, there have been the  Euler-Maclaurin evaluation \cite{Borwein06},
    which also is an analytic continuation of $\zeta$ and the most effective expression in large scale computation.
    We see in computation that the real and image parts of $\zeta(s)$ on the critical line are high-frequency
    oscillation. Even sometimes the two curves are almost tangent and look irregular.
    Corney\cite{Conrey03}(2003) pointed out that
    {\em "It is my belief, RH is a genuinely arithmetic question that
   likely will not succumb to methods of analysis"}, and need more powerful tool.
    Likely proving no zero of the infinite series is hopeless.
    \par
    Next, we turn to $\xi(s)$. Denote $\beta=\sigma-1/2$. Using an asymptotic expansion
     \begin{equation}
      \begin{array}{llll}
    \Gamma(\frac{s}{2})&=\sqrt{2\pi}(\frac{t}{2})^{\beta/2-1/4}e^{-t\pi/4}
     e^{i\phi}(1+O(t^{-1})),\\
    \end{array}
    \end{equation}
    (1.4) and (1.9), there has an important estimate with exponential decay \cite{Borwein06}
    \begin{equation}
    \begin{array}{llll}
    |\xi(s)|& \le C(\frac{t}{2})^{23/12+\beta/6}e^{-t\pi/4}\ln t,\ if\ |\beta|\le 1/2.\\
    \end{array}
    \end{equation}
    Due to the decay $e^{-t\pi/4}$ it is too hard to compute $\xi(s)$ for large $t$.
    Probably this is why there are few work to discuss $\xi$.
    But we can study the geometry property of $\xi(s)$ itself.
    Discussed $\zeta$ and $L$-series, Bombieri\cite{Bombieri}(2000) pointed out that
    {\em "For them we do not have algebraic and geometric models to guide our
    thinking, and entirely new ideas may be needed to study these
    intriguing objects".} He has emphasized the importance to study algebraic and geometric models.
    We begin with it.
    \par
     {\bf Definition 1}. {\em For any fixed $\beta\in (0,1/2]$, $\xi=u+iv$, the sub-interval $I_j=[t_j,t_{j+1}]$ called
   the root-interval, if
   the real part $u(t_j,\beta)=0,u(t_{j+1},\beta)=0$ and $|u(t,\beta)|>0$ inside $I_j$.}
   \par
    {\bf Proposition 1}. {\it For any fixed $\beta\in (0,1/2]$ and in each root-interval
   $I_j=[t_j,t_{j+1}]$, assume that $v(t,\beta)$ has opposite signs at $t_j$ and $t_{j+1}$,
   and $v=0$ at some inner point $t'_j$, then  $\{|u|,|v|/\beta\}$ form local peak-valley structure,
   and norm $||\xi||=|u|+|v/\beta|>0$ in $I_j$, i.e. RH is valid in $I_j$.}
    \par
    We have found analytic and geometric properties of $\xi$, and proved the local peak-valley structure (Theorem 1).
    Because each $t$ must lie in some $I_j$, then $||\xi||>0$ is
    valid for any $t$ (Theorem 2).
    Besides, if RH is valid, based on the equivalence $Re(\frac{\xi'}{\xi})>0$ of
    Lagarias\cite{Lagarias99}(1999), we show that $\xi$ has the peak-valley structure(Theorem 3).
    Therefore both of them are equivalent. \par
     We feel the peak-valley structure may be the geometric model expected by
     Bombieri, which makes the proof of RH get concise and intuitive,
    and many difficulties are avoided, e.g. analyze the summation process of the infinite series
    and prove no zero of it and so on.
    \par

  \section{The $\beta$-symmetry and local peak-valley structure}

   Denote $\tau=it+\beta=s-1/2,\beta=\sigma-1/2$.
  We consider Riemann kernel integral $K(f)$ to define
  \begin{equation}
   \begin{array}{ll}
  \xi(\tau)&=K(f)= \displaystyle \int_1^{\infty}(x^{\tau/2}+x^{-\tau/2})x^{-3/4}f(x)dx=u+iv,\\
  \xi'(\tau)&=K'(f)= \displaystyle \frac{1}{2}\int_1^{\infty}(x^{\tau/2}-x^{-\tau/2})x^{-3/4}\ln xf(x)dx
  =u_{\beta}+iv_{\beta},\\
  \xi''(\tau)&=K''(f)= \displaystyle
  \frac{1}{4}\int_1^{\infty}(x^{\tau/2}+x^{-\tau/2})x^{-3/4}\ln^2xf(x)dx\\
    & =u_{\beta\beta}+iv_{\beta\beta},\\
 \end{array}
   \end{equation}
   which are the alternative high-frequency oscillation. If $\beta=0$, obviously
  \[
  \begin{array}{lll}
  x^{it/2}+x^{-it/2}=2\cos(\frac{t}{2}\ln x),
  \ x^{it/2}-x^{-it/2}=2i\sin(\frac{t}{2}\ln x),
  \end{array}
  \]
  we have the following analytic property.   \par
  {\bf The $\beta$-symmetry}. If  $\beta=0$, then
    \begin{equation}
    \begin{array}{lll}
        v=0,\ \  u_{\beta}=0,\ \  v_{\beta\beta}= 0,\ \ u_{\beta\beta\beta}=0,....
    \end{array}
    \end{equation}
   These properties are essential.
   \par
   {\bf Lemma 1}. {\it  For any $t\in [0,\infty)$ and $\beta\in (0,0.5]$, using the real part $u(t,\beta)$,
    we get the corresponding image part}
    \begin{equation}
    v(t,\beta)= -\int_0^{\beta}u_t(t,r)dr.
    \end{equation}
    \par
    Proof. Using an integral expression
    $v(t,\beta)=v(t,0)+\int_0^{\beta}v_{\beta}(t,r)dr,\ v(t,0)=0,$
    and Cauchy-Riemann condition $v_{\beta}=-u_t$, (2.3) is obtained.
    \par
    {\bf Corollary 1}. $|v(t,\beta)|/\beta$ is uniformly bounded with respect to $\beta\in (0,0.5]$.
    \par
   In the critical strip $S=\{\beta\in (0,0.5],0\le  t<\infty\}$, we define the norm
    \begin{equation}
    ||\xi||=\left\{
    \begin{array}{lll}
    |u|+|v|/\beta, & if \ \beta \in (0,1/2],& t\in [0,\infty),\\
    |u(t,0)|+|u_t(t,0)|,& if\ \beta\rightarrow +0,& t\in [0,\infty),\\
    \end{array}
    \right.
    \end{equation}
    where three conditions of norm are satisfied.
    The advantage is that $|u|$ and $|v|/\beta$ are of the same order and $||\xi||$
    is stable with respect to $\beta>0$.
     \par
 Finally we want to explain the local peak-valley structure by the curve figures with
 $\beta= 0.1, 0.3$ and $0.5$ in Fig.1.1-4, where we have used a changing scale $M=8(t/2)^{23/12+\beta/6}e^{-t\pi/4}$ and $(u/M,v/M)$ when drawing these curves.
 Fig.1.2 shows that $u(t_2,\beta)=0,u(t_3,\beta)=0$ at two end-points of root-interval $I_2=[t_2,t_3]$, and $u(t,\beta)>0$ inside $I_2$. We also see that $v(t_2,\beta)<0$,$v(t_3,\beta)>0$, and $v(t'_2,\beta)=0$ at some inner point $t'_2$ of  $I_2$. In Fig.1.3, there has a peak for $|u|\ge 0$, while  has a valley for $|v|/\beta\ge 0$, then it forms a peak-valley structure of $\{|u|,|v|/\beta\}$ in $I_2$. Fig.1.4 exhibits the low bound $\min_{t\in I_2}(|u|+|v|/\beta)/M\ge 0.0876$, i.e. RH is valid in $I_2$.
\begin{figure}
	\centering 
	\includegraphics[height=168pt]{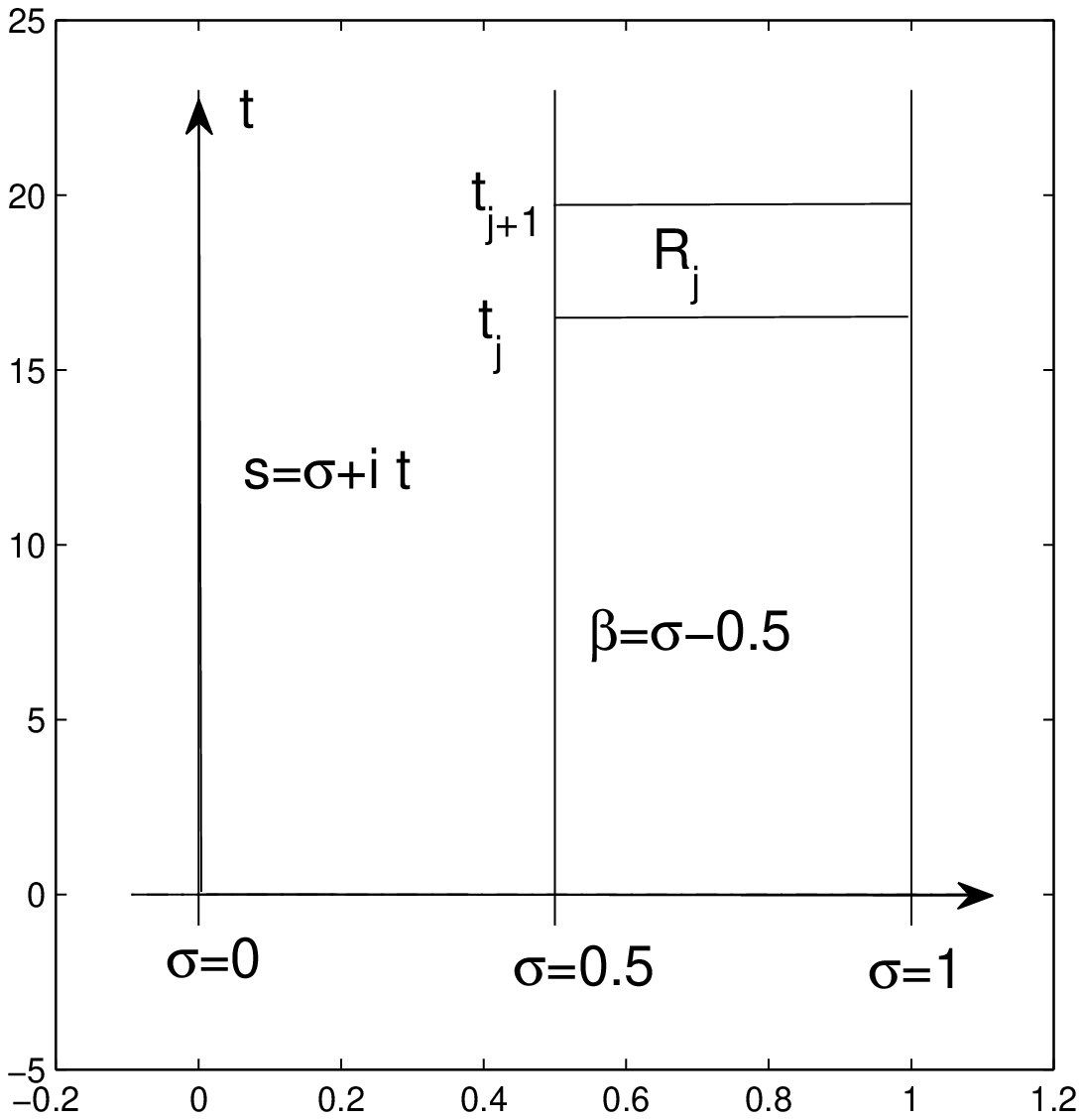}
	\includegraphics[height=168pt]{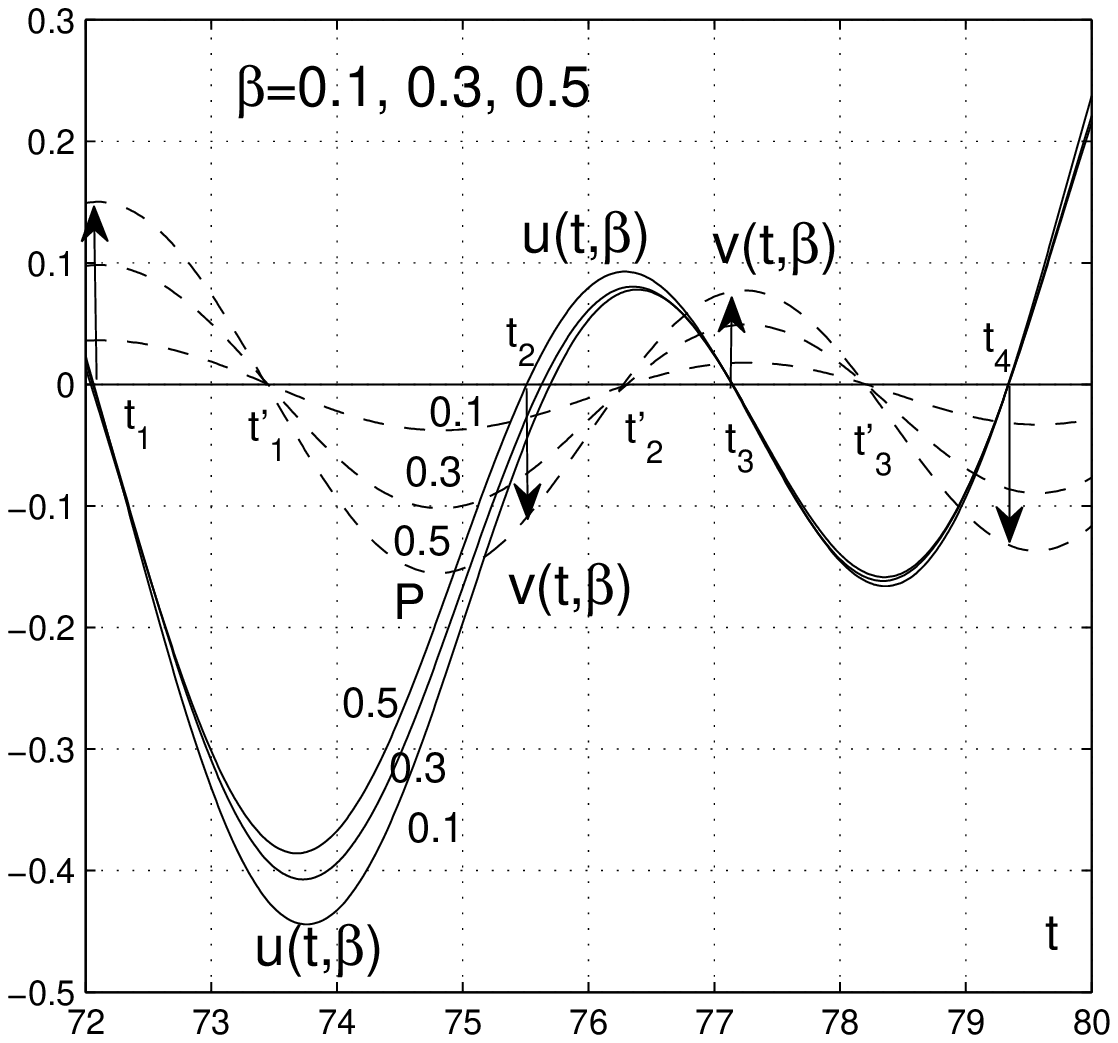}
	\includegraphics[height=160pt]{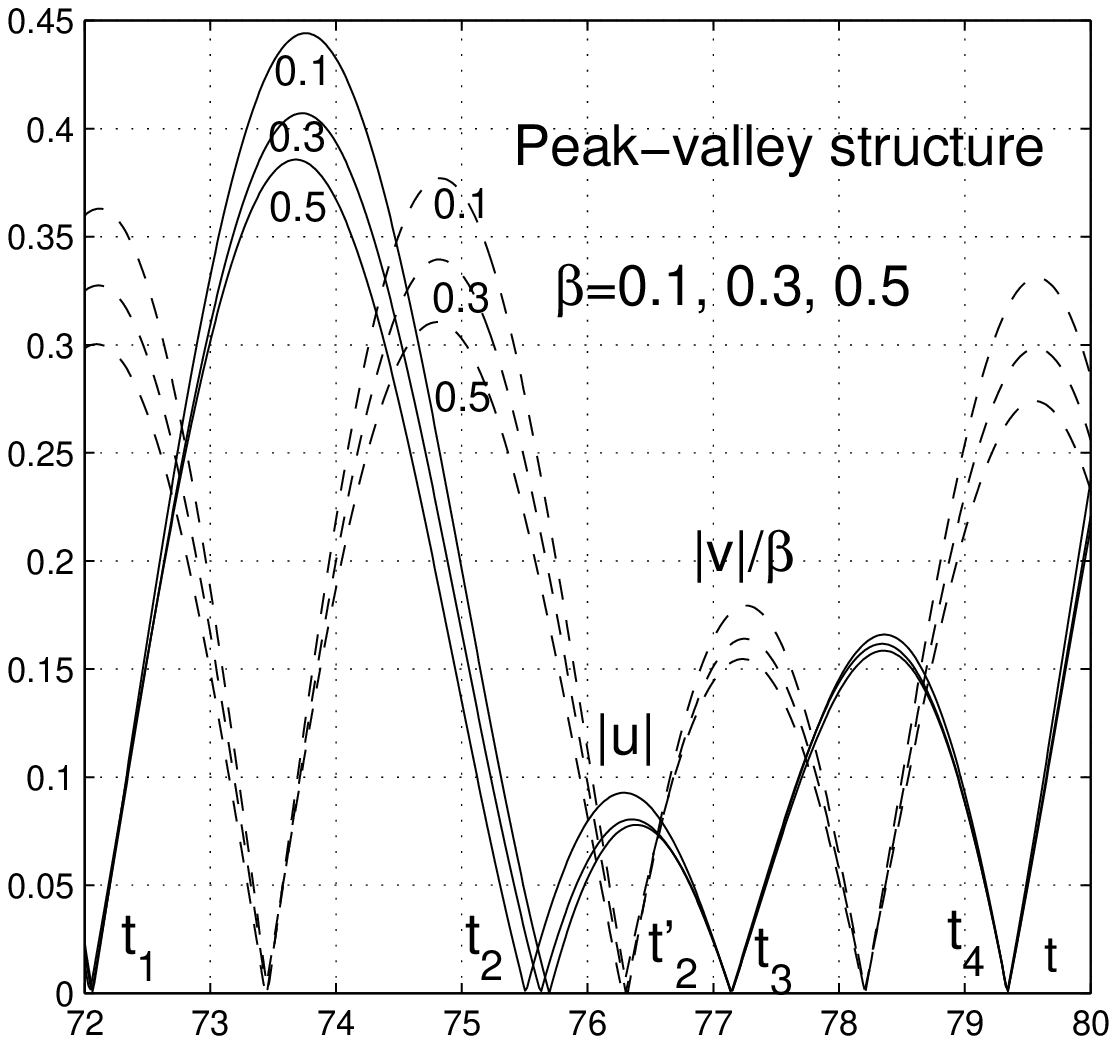}
	\includegraphics[height=160pt]{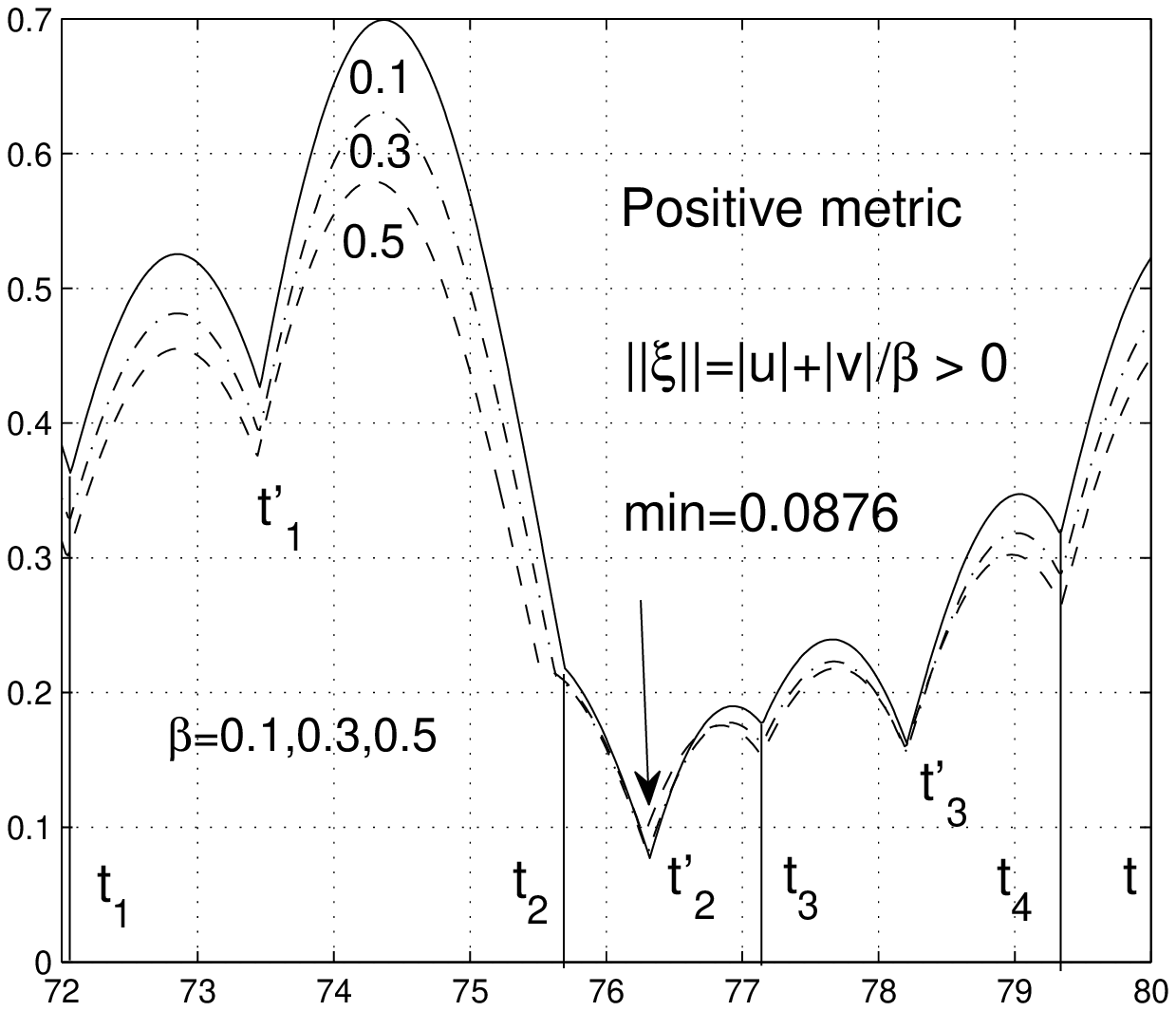}
	\caption{Fig. 1.1-4. Take $\beta=0.1,0.3,0.5$, curves $\{u,v\}$,$\{|u|,|v|/\beta\}$ and $|u|+|v|/\beta>0$.}
\end{figure}

\section{Local geometric proof of RH}
 We shall regard $\{u(t,\beta),v(t,\beta)\}$ as a continuous changing process
   from $\beta=+0$ to $\beta=0.5$.  For any fixed $\beta\in (0,0.5]$ the real part $u(t,\beta)$
   is an irregular high-frequency oscillation, and its zeros $t_j$
   (also depend on $\beta$) form an irregular infinite sequence
   \[
   \begin{array}{lll}
   ...<t_{j-1}<t_j<t_{j+1}<t_{j+2}<.....
   \end{array}
   \]
   We shall take them as the base in studying peak-valley
   structure. We prove
   \par
   {\bf Theorem 1(local peak-valley structure)}. {\it For any fixed $\beta\in (0,1/2]$ and in each root-interval
   $I_j=[t_j,t_{j+1}]$, then the curves $\{|u|,|v|/\beta\}$ form a local peak-valley structure,
   and $||\xi||$ has the positive low bound independent of $t\in  I_j$,}
    \begin{equation}
     \min_{t\in I_j}(|u|+|v|/\beta)=\mu(t_j,\beta)>0,\ \ \ \beta\in [+0,0.5].
    \end{equation}
    \par
    Proof. {\bf 1. Single peak case}.\par
    From Fig.1.2 we have seen the following general property.
    \par
    {\bf Geometric property of single peak}. {\em For any $\beta\ge 0$, there are $u_t>0$
    from negative peak to positive one, and $u_t<0$ from positive peak
    to negative one.}
    \par
     Below it is enough to discuss $u>0$ inside the root-interval $I_j$. For any fixed
    $\beta>0$, using Lemma 1, we discuss two cases as follows.\par
   1). As $u_t>0$ near the left node $t_j$,we have
   \begin{equation}
   \left\{
   \begin{array}{ll}
   v(t_j,\beta)/\beta&=- \displaystyle \frac{1}{\beta}\int_0^{\beta}u_t(t_j,r)dr<0,\\
   \lim\limits_{\beta\rightarrow +0}v(t_j,\beta)/\beta&=-u_t(t_j,0)<0.
   \end{array}
   \right.
   \end{equation}
     \par
   2). As $u_t<0$ near the right node $t_{j+1}$, similarly
  \begin{equation}
  \left\{
   \begin{array}{ll}
    v(t_{j+1},\beta)/\beta&=   \displaystyle  -\frac{1}{\beta}\int_0^{\beta}u_t(t_{j+1},r)dr>0,\\
   \lim\limits_{\beta\rightarrow +0}v(t_{j+1},\beta)/\beta&=-u_t(t_{j+1},0)>0.\\
   \end{array}
   \right.
   \end{equation}
   which are valid and numerically stable for $\beta\in [+0,1/2]$.
   \par
   Because $v(t,\beta)$ has opposite signs at two end-points in $I_j$,
   there certainly exists an inner point $t'_j=t'_j(\beta)$ such that $v(t'_j,\beta)=0$.
   Clearly in $I_j$, $|u|$ is a peak curve and $|v(t,\beta)|/\beta$ is a valley curve,
   thus $\{|u|,|v|/\beta\}$ form a local peak-valley structure.
   We define a continuous function with respect to $(t,\beta)$
   \[
   \begin{array}{lll}
    \phi(t,\beta)=|u(t,\beta)|+|v(t,\beta)|/\beta,\ \ \ \beta\in [+0,0.5],\ t\in I_j=[t_j,t_{j+1}],
    \end{array}
    \]
   which certainly has a positive low bound independent of $t\in I_j$,
    \begin{equation}
     \min_{t\in I_j}\phi(t,\beta)=\mu(t_j,\beta)>0,\ \ \beta\in [+0,0.5].
    \end{equation}
    This is a fine local geometric analysis.
    \par
   {\bf 2. Multiple peak case}. Although $||\xi||>0$ is still valid, we shall prove that
   multiple peak case does not appear.
   \par
     Assume that $u(t,\beta)>0$ inside $I_j=[t_j,t_{j+1}]$, and there has odd number of positive
     extreme values $u(t_{jp},\beta)=a_{jp}>0$ at the inner points $t_{jp}$,
    $p=1,2,...,2k+1$. We consider the minimum value $u(t_{ji},\beta)$ at some point
    $t'=t_{ji}$, where $u>0, u_t=0$ and $u_{tt}>0$, i.e., $u$ is convex toward $t$-axis.
     By Cauchy-Riemann condition $v_t=u_{\beta}$ and analytic property $u_{\beta}(t,0)=0$, we get
    \begin{equation}
     v_t(t',\beta)=v_t(t',0)+\int_0^{\beta}v_{t\beta}(t',r)dr= -\int_0^{\beta}u_{tt}(t',r)dr<0.
     \end{equation}
    On the other hand, because $u>0$ inside $I_j$, then RH is
    locally valid. Using the equivalence of Lagarias \cite{Lagarias99}(1999), we have
    \begin{equation}
    Re(\frac{\xi'}{\xi})=Re(\frac{\xi'\bar{\xi}}{|\xi|^2})
    =(uu_{\beta}+vv_{\beta})/|\xi|^2>0,\ \ \beta>0.
    \end{equation}
    Using Cauchy-Riemann conditions $u_{\beta}=v_t,u_t=-v_{\beta}$, it derives
    \begin{equation}
     \psi(t)=uv_t-vu_t>0,\ \ \beta>0.
    \end{equation}
    But now, $u>0, u_t=0, v_t<0$ at $t=t'$, which lead to contradiction $\psi=uv_t<0$.
    \par
    Should point out that in multiple peak case for fixed $\beta>0$, the real part $u\ge  a_{jl}>0$ in the subinterval
    $I^*=[t_{j1},t_{j,2k+1}]$, while in remaining subintervals $[t_j,t_{j1}]$ and $[t_{j,2k+1},t_{j+1}]$,
     we can still get the positive estimates (3.2) and (3.3),
    therefore $||\xi||>0$ in the root-interval $I_j$.
    But the multiple peak case will imply the following danger: When $\beta>0$ is further increased,
    the minimum extreme value $u=a_{ji}>0$ possibly gradually is close to $t$-axis toward its convex direction,
    as seen in the case 3, we can not deny the possibility to contact with $t$-axis,
    which will lead more difficulty.
    Fortunately, we have proved that the multiple peak case does not appear by using the equivalence theorem
    of Lagarias, and extricated oneself from the difficult position.
    \par
    {\bf 3. The case of two zeros to be very close to each other on the critical line}.
    We shall prove that the root-interval will be enlarged for $\beta>0$ rather than
    decreased, and the peak-valley structure is valid. This conclusion is also valid in the case of double root,
    although no double root is found in computation up to now.
    \par
    Assume that $u(t,0)$ on the critical line has a solitary small
    interval $I^0=[t^0_j,t^0_{j+1}]$, $u(t^0_j,0)=0$ and $u(t^0_{j+1},0)=0$.
    There is an extreme point $t'\in I^0_j$ such that
    $u( t',0)=\epsilon >0$ and $u_t=0$. We consider a larger interval $I\supset I^0_j$,
    in which $u_{tt}<0$ and $u(t,0)$ is convex upwards, then $u_t>0$ for $t<t'$ and $u_t<0$ for
    $t<t'$. We say $\epsilon>0$ very small if
    $\epsilon/|u_{tt}|<<1$. An artificial example is shown in Fig.2.
    \par
    Now consider small $\beta>0$, in the larger interval $I$, we have
   \begin{equation}
   \begin{array}{lll}
   u(t,\beta)&= \displaystyle  u(t,0)+u_{\beta}(t,0)\beta-\int_0^{\beta}u_{\beta\beta}(t,r)(r-\beta)dr&as\ u_{\beta}(t,0)=0\\
         &=  \displaystyle u(t,0)+d,\ d=-\int_0^{\beta}u_{tt}(t,r)(\beta-r)dr>0 &as\ u_{\beta\beta}=-u_{tt}\\
   v(t,\beta)&=  \displaystyle -\int_0^{\beta}u_t(t,r)dr,\ \ \ \ see\ (2.3)   &as\ v_{\beta\beta}(t,0)=0\\
          &=   \displaystyle -u_t(t,0)\beta+\frac{1}{2}\int_0^{\beta}u_{ttt}(t,r)(r-\beta)^2dr,& as\ v_{\beta\beta\beta}=u_{ttt},
    \end{array}
   \end{equation}
   which can be summarized in the following form ($\epsilon=0$ is admissible)
   \begin{equation}
   \left\{
   \begin{array}{lll}
   t=t'& u(t',\beta)=\epsilon+d>0,       & v(t',\beta)=O(\beta^3),\ \ as\ u_t(t',0)=0,\\
   t<t'&u(t,\beta)=u(t,0)+d,           &v(t,\beta)= \displaystyle  -\int_0^{\beta}u_{t}(t,r)dr<0,\\
   t>t'&u(t,\beta)=u(t,0)+d,  & v(t,\beta)= \displaystyle  -\int_0^{\beta}u_{t}(t,r)dr>0,\\
    \end{array}
   \right.
   \end{equation}
  From this we see that $u(t,\beta)$ has removed $u(t,0)$ in parallel
  by a distance $d>0$ toward the direction of its convexity.
  Due to $u(t,0)<0$ outside $I^0$, there are certainly a left node
  $t_j=t_j(\beta)$ and a right node $t_{j+1}=t_{j+1}(\beta)$ such
  that $u(t,\beta)=0$. Inside the enlarged interval $I_j=[t_j,t_{j+1}]$,
  it forms a positive peak curve for $u(t,\beta)>0$. Besides, $v(t,\beta)$
  has opposite signs at two endpoints of $I_j$, and there certainly exists some
  inner point $t''$ such that $v(t'',\beta)=0$, i.e., $|v|/\beta$ is
  a valley curve. Therefore it forms a peak-valley
  structure for $\{|u|,|v|/\beta\}$ in $I_j$.

  \par
  It should be pointed out that if  $u_{tt}<0$ and decreasing $u(t',0)< 0$, which belongs
  to the multiple peak case. This case does not appear as mentioned above.
  \par
  Finally by summarizing three cases,  Theorem 1 is proved.
  \par

  \begin{center}
  \includegraphics[height=168pt]{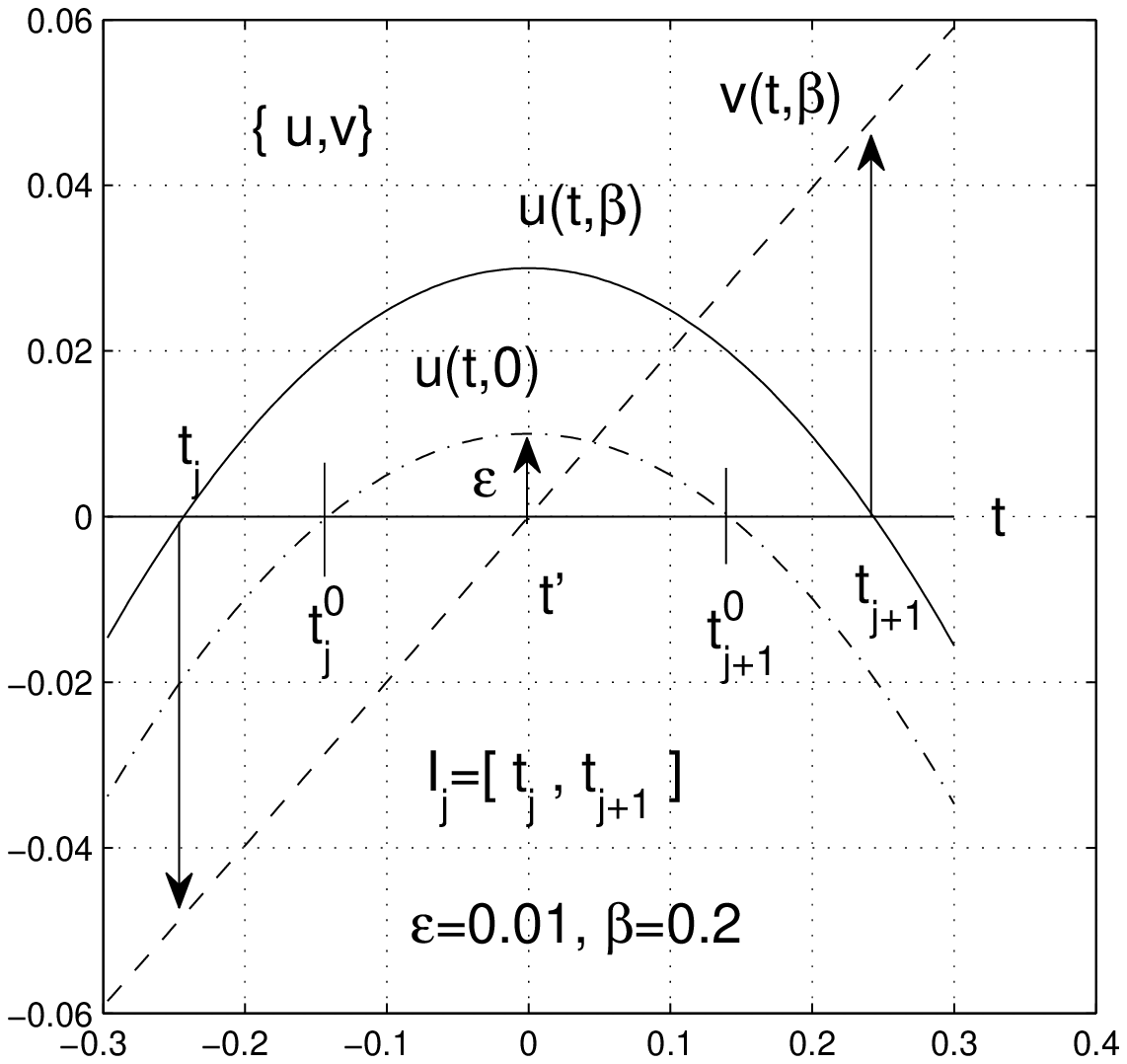}
  \includegraphics[height=168pt]{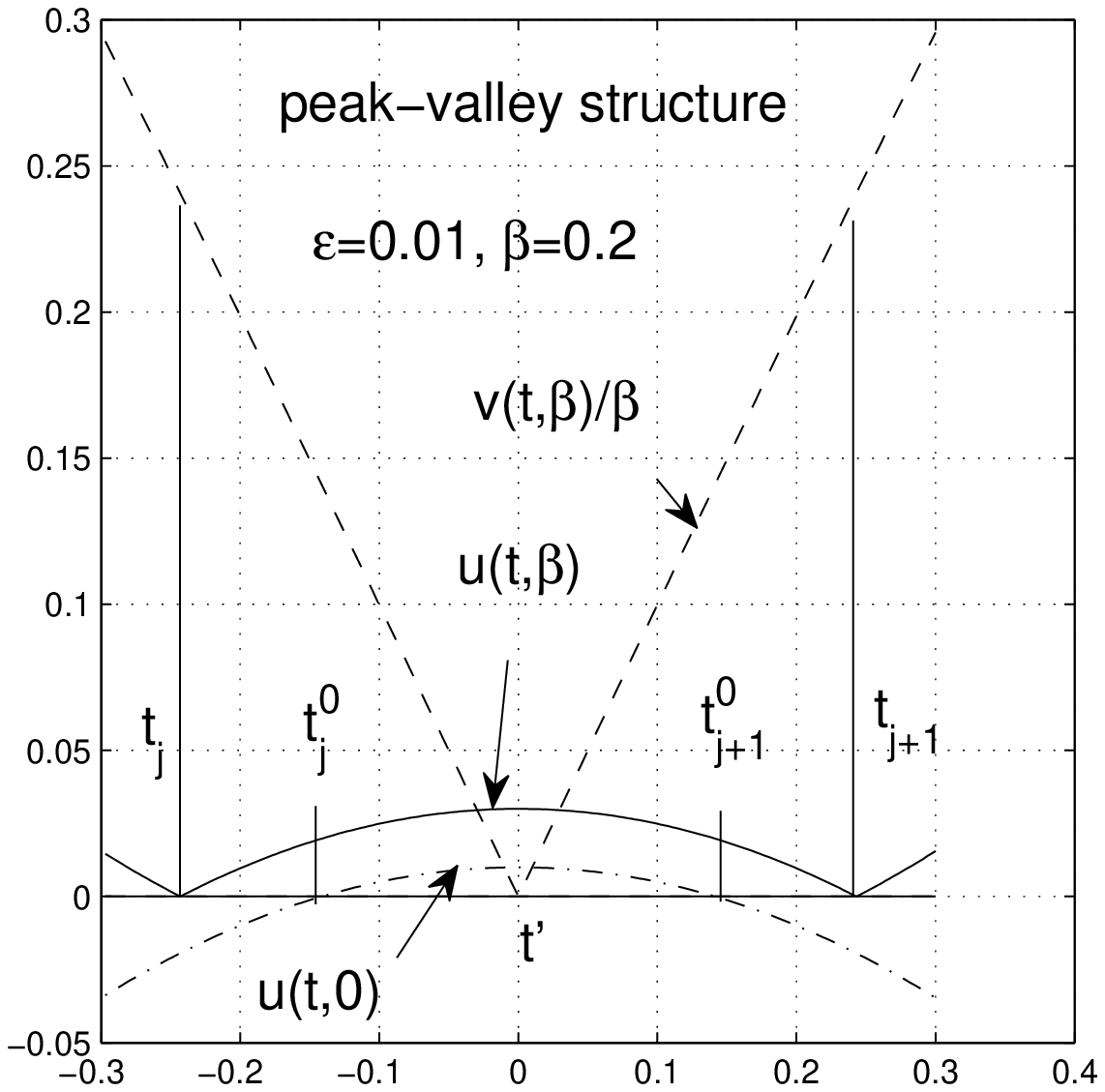}
    \begin{center}
   Fig.2. Initial value $u(t',0)=\cos(t)-1+\epsilon$, peak-valley structure for $\beta=0.2$\\
   \end{center}
   \end{center}

    We have constructed an example $u(t,0)=\cos(t)-1+\epsilon, \epsilon=0.01, \beta=0.2$ in Fig.2,
   which has a local peak-valley structure in $I_j=[t_j,t_{j+1}]$, and it is also valid for $\epsilon=0$
   which is double root.
   Besides, we have seen in Fig.1.2 that when increasing $\beta$, the corresponding smaller interval
   $(t_2,t_3)$ is enlarged, while another neighbor larger interval $(t_1,t_2)$ will be decreased.
   \par

   {\bf Theorem 2}. {\it RH is valid for any $(\beta,t)\in (0,0.5]\times [0,\infty)$.}
    \par
    Proof. Actually, for any fixed $\beta\in (0,0.5]$, an infinite sequence can be formed for the zeros
    $\{t_j(\beta)\}$ of $u(t,\beta)$, while any $t\in [0,\infty)$ must lie in some $I_j$ such that $||\xi||\ge \mu(t_j,\beta)>0$.
   The theorem is proved.\par
    Recall that the equivalence $Re(\frac{\xi'}{\xi})>0$ of Lagarias in \cite{Lagarias99}(1999),
    which is a unique equivalent theorem with $\xi'$ up to now. We will prove the following.
   \par
   {\bf Theorem 3}. {\it The peak-valley structure and RH are  equivalent.}
   \par
   { Proof}. Assume that RH is valid and $u>0$ inside root-interval $I_j=[t_j,t_{j+1}]$(similarly for $u<0$).
    By (3.7), $\psi=uv_t-vu_t>0$, for $\beta>0$. We have the following facts.\par
   At the left node $t_j$, $u=0, u_t>0$(geometric property) and $\psi=-vu_t>0$, then $v<0$;\par
   At the right node $t_{j+1}$, $u=0, u_t<0$ and $\psi=-vu_t>0$, then $v>0$. \par
    Therefore $v$ has opposite signs at two end-points, there certainly exists
    an inner point $t'_j \in I_j$ such that $v=0$, which implies local peak-valley
    structure. Theorem 3 is proved.
    \par
    From the view-point of complex analysis, RH requires $|\xi|>0$,
    while from the view-point of geometry, the peak-valley structure
    requires stronger norm $||\xi||>0$. Both of them are equivalent.
    However, the local geometry property is of extreme importance,
    because proving the peak-valley structure is concise and intuitive.\par
     {\bf Remark 1}. {\em In the proof of Theorem 1 we have seen  that
     Riemann integral $\xi=K(f)$  has $\beta$-symmetry, which is independent
     of the speciality of $f$. Therefore we guess that for the very wide class of
     the fast decay function $f$, RH is still valid for $K(f)$.}
       We have two examples. For $t\le 110$, there are larger low bounds
      $(|u_{\beta}|/\beta+|v_{\beta}|)/M\ge 0.20$ and $||\xi''||/M\ge 0.28$.
     \par
      Haglund\cite{Haglund} has discussed $\xi$ and other functions with numerical experiments,
      and proposed a conjecture: If function $F_N$ has monotonic zeros,
      then which implies RH. Sarnak\cite{sarnak} has analyzed the Grand RH of L-function,
      which is more difficult.


\begin{center}
   {\bf References}
   \end{center}
 \begin{enumerate}
   {\small

     \bibitem{Bombieri} E. Bombieri. Problems of the Millennium: The
     Riemann Hypothesis. AMS. 107-124(2000) \\
     http://www.claymath.org.

   \bibitem{Borwein06} P. Borwein, S. Choi, B. Rooney, A.
    Weirathmuller. The Riemann Hypothesis. Springer (2006)


   \bibitem{Conrey03} J. Conrey. The Riemann Hypothesis.
   Notices of The AMS. 341-353(2003)


    \bibitem{Edwards} H. Edwards. Riemann's Zeta function. Mineola:
    Dover Publication,Inc.(2001)


    \bibitem{Haglund} J. Huglund. Some conjectures on the zeros of
    approximates to the Riemann zeta-function and incomplete gamma
    functions. Cent. Eur. J. Math., {\bf 9:2},302-318(2011)

   \bibitem{Lagarias99} J. Lagarias. On a positivity property of the Riemann
    $\zeta$-function. Acta Arith. {\bf 89:3}, 213-234(1999)


    \bibitem{Lune83} J. Lune, H. Riele. On the zeros of the Riemann
     Zeta function in the critical strip. Part 3. Math. Comp., {\bf
     41}, 759-767(1983)

     \bibitem{Lune86} J. Lune, H. Riele, D.Winter. On the zeros of the Riemann
     Zeta function in the critical strip. Part 4. Math. Comp., {\bf 46}, 667-681(1986)

     \bibitem{sarnak} P. Sarnak. Problems of the Millennium: The
      Riemann Hypothesis (2004).\\      
       http://www.claymath.org.


    \bibitem{Smale00} S. Smale. Mathematical problems for next century.
      Math. Intelligencer. {\bf 20:2}, 7-15(1998) 

      }

     \end{enumerate}
\end{document}